\documentclass[a4paper]{article}

\author{Marco Leoni}
\date{}
\begin{document}
\title{Chow groups of weighted hypersurfaces.}

\maketitle

{\bf Abstract }| We extend a result of to Esnault-Levine-Viehweg concerning
  the Chow groups of hypersurfaces in projective  space to those in
  weighted projective spaces. 

\newtheorem{theo}{Theorem}[section]
\newtheorem{lem}{Lemma}[section]
\newtheorem{ex}{Example}[section]
\newtheorem{defin}{Definition}[section]
\newtheorem{rem}{Remark}[section]
\newtheorem{cor}{Corollary}[section]
\newtheorem{cl}{Claim}[section]

\section{Introduction}
The purpose of this paper is to generalize to the case of weighted 
projective spaces over an algebraically closed field ${\bf K}$ the 
following result from [ELV]:
\begin{theo}{\rm [ELV], Th. 4.6.}
Let $X\subset {\bf P}^n$ be a hypersurface of degree $d\geq 3$ and let
$s\leq n-1$ be an integer such that:
$$\left(
\begin{array}{c}
s+d\\
s+1\end{array}\right) \leq n.$$
Then ${\rm CH}_s (X)\otimes {\bf Q}={\bf Q}$.
\end{theo}

Let $Q=(q_0 ,\dots ,q_n )\in {\bf N}^{n+1}$. Let $\mu_a :=
\{z\in {\bf K}\vert z^a =1\}$ and set $\mu :=\prod_{i=1}^n \mu_{q_i}$.
The weighted projective space 
${\bf P}(Q)$ can be realized either as the quotient ${\bf P}^n /\mu$
(with the action defined by componentwise multiplication) 
or as the quotient ${\bf K}^{n+1}/\bf{K}^*$, the action being defined by 
$t\cdot (x_0 ,\dots x_n ):=(t^{q_0}x_0 ,\dots ,t^{q_n}x_n )$.
The map $\varphi_Q:[t_0 :\dots 
:t_n]_\mu \rightarrow [t^{q_0}_0,\dots ,t^{q_n}_n ]_{{\bf K}^*}$ 
gives the isomorphism between the two representations. 
One deduces from this that 
there is a one to one correspondence between hypersurfaces $X:=\{ f=0\}$
of ${\bf K}^{n+1}/{\bf K}^*$ and those of ${\bf P}^{n}/\mu$ defined by the
zeroes of the polynomial $f'([t_0 :\dots ;t_n ]_\mu ):=f([t_0^{q_0},\dots,
t^{q_n}_n ]_{{\bf K}^*} )$.
If $f'$ is smooth, the hypersurface $\{ f=0 \}$ is seen to be 
{\it quasismooth}: the cone ${\cal C}_X :=\{ x\in {\bf K}^{n+1}\vert f(x)=0\}$
has one singularity in the origin.
\section{The Main Result.}
\begin{theo}
For a smooth irreducible weighted hypersurface $X'$ of degree 
$d\geq 3$ in ${\bf P}^n$ and $\forall l\in {\bf N}$ such that:
$$\left( \begin{array}{c}
d+l\\
1+l\end{array}\right)\leq\sum_{j=0}^n q_j -1$$
one has:
$${\rm CH}_l (X)\otimes {\bf Q}={\bf Q}$$
where $X=X'/\mu$. 
\end{theo}

{\bf Proof:}

Let:
$$N:=(\sum_{j=0}^n q_j )-1$$
$$N_r :=\left\{ \begin{array}{cc}
0 &  \, r=-1\\
\sum_{j=0}^r q_j &\forall r=0,\dots, n
\end{array}\right.$$
Remark in particular that $N_n =N+1$, and that $N_r-N_{r-1}=q_r$
$\forall r=0,\dots ,n$.
Define a rational map:
$$\sigma_Q :{\bf P}^N \rightarrow{\bf P}(Q)$$
by:
$$\left( \sigma_Q ([t_0 :\dots :t_N ]\right)_r :=\prod_{j=N_{r -1}}^{N_r 
-1}t_j \,\, \forall r=0,\dots ,n.$$

Set: 
$${\cal J}_Q :=\{ (j_0 ,\dots ,j_n )\in{\bf N}^{n+1}:N_{r-1}\leq j_r \leq
N_r -1 \, \forall r=0,\dots ,n \}$$
and consider $\forall J\in {\cal J}_Q$, the subvarieties:
$$Z_J :=\left\{ t\in{\bf P}^N :t_{j_0}=\dots =t_{j_n}=0\right \}$$
$$Z_Q :=\cup_{J\in{\cal J}_Q }Z_J$$
It is clear that $\sigma_Q$ is only defined on ${\bf P}^N -Z_Q$.
This map is well-defined on ${\bf P}^N -Z_Q$: indeed, if one
considers $lt_j$ instead of $t_j$ for a nonzero $l$, one has:
$$\prod_{j=N_{r-1}}^{N_r -1}lt_j =l^{q_r}\prod_{j=N_{r-1}}^{N_r -1}t_j$$
so that modulo the weighted action of ${\bf K}^*$ these two quantities 
coincide. 

Also, $\sigma_Q$ is onto, since if $x\in{\bf P}(Q)$ and
$(x_0 ,\dots ,x_n )$ is a representative in ${\bf K}^{n+1^*}$, one may
choose, $\forall r=0, \dots ,n$, some $q_r -1$ variables freely and 
the last one such that $x_r =\prod_{j=N_{r-1}}^{N_r -1}t_j $. So:
$$\forall x\in {\bf P}(Q), \dim (\sigma_Q ^{-1}(x))=\sum_{r=0}^n (q_r -1)=
N-n$$

Let $X\subset {\bf P}(Q)$ be a weighted homogeneous hypersurface of 
$Q$-degree $d\geq 3$. If $X$ is defined by the weighted homogeneous polynomial
$f=f(x_0 ,\dots ,x_n )$, we define $\tilde{X}$ in ${\bf P}^N$ by the 
polynomial $\tilde{f}=\tilde{f}(t_0 :\dots :t_N )$, of the same degree, 
obtained by replacing $x_k$ by $\prod_{j=N_{k-1}}^{N_k -1}t_j$. 
The map $\sigma_Q$ induces a rational map:
$$\sigma_Q :\tilde{X}\rightarrow X.$$

Let $R$ be the plane in ${\bf P}^N$ defined by the equations:
$$t_{N_{r-1}}=\dots =t_{N_r -1}\quad \forall r=0,\dots ,n$$
The number of equations which define it is:
$$\sum_{r=0}^n (N_{r}-1-N_{r-1})=\sum_{r=0}^n q_r -(n+1)=N-n$$
Let $S:=R\cap \tilde{X}$. Then this linear space has dimension $n$
and has, by construction, the fundamental property that
$S\cap Z_Q =\emptyset$:
$$\forall t\in Z_Q ,\forall r\vert 0\leq r\leq n,\exists i\mbox{ such that }
N_{r-1}\leq i \leq N_r -1\mbox{ for which } t_i =0$$
But then in $S$, $t_{N_{r-1}}=0$ also and all the other $t_j$ with $j$
in the $r$th string are also zero. This for every $r$.

Let 
$$u:{\rm Bl}_{Z_Q}(\tilde{X})\rightarrow \tilde{X}$$
be the blow-up along $Z_Q$ turning $\sigma_Q$ into a morphism:
$$\begin{array}{ccc}
{\rm Bl}_{Z_Q}(\tilde{X})&\stackrel{\hat{\sigma}_Q}{\rightarrow}&X\\
\downarrow u & &\vert\vert \\
\tilde{X}&\stackrel{\sigma_Q}{\rightarrow}& X
\end{array}$$
Let:
$$l_0 := \max_{l\in {\bf N}} \left\{ \left( \begin{array}{c}
l+d\\
l+1\end{array}\right) \leq N\right\}\cap \{ l\in {\bf N}\vert l\leq n\} $$

We know from [ELV], Theorem 4.6., that if $s\leq l_0$, then: 
$${\rm CH}_s (\tilde{X})\otimes{\bf Q}={\bf Q}$$

So let's take $\gamma\in {\rm CH}_s (X)\otimes {\bf Q}$ where $s\leq l_0$.
Set $\tilde{\gamma}:=\tau^{-1}(\gamma )$, being $\tau :=\sigma_Q \vert_S$.

Certainly $\tilde{\gamma}$ is an $s$-cycle on $\tilde{X}$ which 
is supported on $S$.
Therefore there is some $a\in{\bf Q}$ and a $\Gamma\in {\rm Gr}(s+1)$
such that:
$$\tilde{\gamma}\sim_{\tilde{X}}[\Gamma\cap\tilde{X}]=a\Gamma\cdot \tilde{X}$$
Since ${\rm CH}_s ({\bf P}^N )\otimes {\bf Q}
={\bf Q}$, one can eventually replace $\Gamma$ by another 
$(s+1)$-plane which is transversal to $Z_Q$. 
Therefore we may assume that the proper transform of
$\Gamma$ under the blow-up along $Z_Q$, which I denote by 
$\hat{\Gamma}$, is isomorphic to $\Gamma$ itself.
Certainly $\hat{\tilde{\gamma}}\simeq \tilde{\gamma}$ because $\tilde{\gamma}
\cap Z_Q =\emptyset$.
Therefore we deduce:
$$\hat{\tilde{\gamma}}\sim_{{\rm Bl}_{Z_Q}(\tilde{X})}b\hat{\Gamma}\cdot
{\rm Bl}_{Z_Q}(\tilde{X})$$
Since $X'$ is smooth, and since $\mu$ is a finite group, by [FU], Ex 11.4.7.,
we have a ``moving lemma'' on $X=X'/\mu$. Therefore we can move $\gamma$
inside $X$ in such a way that that it is not in the ramification locus of
$\hat{\sigma}_Q$. Hence $\hat{\sigma}_Q$ is finite of a certain nonzero degree,
say $e$. 
%\footnote{Indeed:
%$$\tau^{-1}(x_0 ,\dots ,x_n )=\{ [t_0 :\dots :t_0 :\dots :t_n :
%\dots :t_n ]\in {\bf P}^N \vert t^{q_j}_j =x_j \, \forall j\}$$
%so $\tau$ is finite.}

So we deduce:
$$\hat{\sigma}_{Q*}\hat{\tilde{\gamma}}=e\gamma$$ 
while:
$$\hat{\sigma}_{Q*} (\hat{\Gamma '}\cdot {\rm Bl}_{Z_Q}(\tilde{X}))=
eH_{s+1}\cdot eX$$
being $H_{s+1}$ the generator of ${\rm CH}_{s+1}({\bf P}(Q))\otimes
{\bf Q}={\bf Q}$.

Therefore $\gamma\sim_{X}be^2 H_{s+1}\cdot X=tH_s$, with $H_s$ generator of
${\bf P}(Q)\otimes{\bf Q}={\bf Q}$. This shows ${\rm CH}_s(X)\otimes
{\bf Q}={\bf Q}$ $\forall s\leq l_0$.
\begin{flushright}
QED
\end{flushright}
\begin{rem}
\par\rm
In the preceding proof a moving Lemma is used; for this reason
  $X'$ should have at most quotient singularities. One can probably
  avoid this as to arrive at the true generalization of the result
  in [ELV] valid irrespective of the singularities.

Essentially the same method also works for complete intersections
  so that appropriate analogues of [ELV] Prop. 3.5 and
  Thm.~4.6. hold. In view of technical complications we preferred to
  state and give the proof for hypersurfaces onl

\end{rem}

{\vspace{0.5cm}

{\bf References:}
\newline
[ELV] H.Esnault- M.Levine- E.Viehweg ``{\it Chow groups of projective 
varieties of very small degree}'' Duke Math. Journal 87 n.1 (1997) 29-58.
\newline
[FU] W.Fulton ``{\it Intersection Theory}'', Springer-Verlag, Berlin, 1984.

\end{document}